\documentclass{birkjour}

\usepackage{amscd}
\usepackage{amssymb}
\usepackage[all,cmtip]{xy}

\newtheorem{thm}{Theorem}[section]
 \newtheorem{cor}[thm]{Corollary}
 \newtheorem{lem}[thm]{Lemma}
 \newtheorem{prop}[thm]{Proposition}
 \theoremstyle{definition}
 
 \newtheorem{rem}[thm]{Remark}
 \newtheorem{nada}[thm]{}
 
  \theoremstyle{remark}
 
 \numberwithin{equation}{section}

\begin{document}

\title[Linear-type maps]{LINEAR-TYPE MAPS BETWEEN AFFINE AND PROJECTIVE SPACES}

\author[Juan B. Sancho]{Juan B. Sancho de Salas}

\address{Dpto. de Matem\'{a}ticas, Univ. de Extremadura, E-06071 Badajoz (SPAIN)}

\email{jsancho@unex.es}



\date{}

\begin{abstract}
We determine those maps between affine or projective spaces that are {\it linear} in the abstract sense of transforming collinear points into collinear points and whose restriction to any line is constant or injective. Our results are extensions of the classical Fundamental Theorems of Affine and Projective Geometries.
\end{abstract}

\maketitle

\section*{Introduction}

Classically, the study of affine and projective geometries focused its attention on isomorphisms (affinities and projectivities). This seems very restrictive from the modern categorical point of view: What is the category where we must place affine and projective spaces and other incidence geometries? Texts on affine and projective geometry usually ignore such a question. For example, the basic fact that any affine space admits a natural embedding into a projective space is presented (see \cite{Bennett} or \cite{Berger}) without proposing a definition of embedding.

To set the appropriate category, let us observe that the structure of affine and projective spaces is determined by their lattices of linear subspaces. They are examples of {\it closure spaces} satisfying MacLane's exchange axiom and a finitary condition. For the sake of brevity, following \cite{FF} and \cite{Gratzer}, such objects are simply named {\it geometries}. These objects have a natural notion of morphism: A map $\varphi\colon X\to Y$ is a morphism if $\varphi^{-1}S$ is a subspace of $X$ for any subspace $S$ of $Y$. In the most common case in this paper (geometries generated by lines), a morphism of geometries $\varphi\colon X\to Y$ is just a map transforming collinear points into collinear points and such that it is injective or constant on any line.

Our purpose is to algebraically characterize morphisms of geometries between affine or projective spaces, in terms of affine or homogeneous coordinates. Our results are extensions of the classical Fundamental Theorems of Affine and Projective Geometries.

 In the case of projective spaces, the question is resolved by a result of Faure--Fr\"{o}licher \cite{FF1} and Havlicek \cite{Havlicek}, which is a natural generalization of the classical Fundamental Theorem of Projective Geometry to non necessarily bijective maps:\medskip

\noindent\textbf{Fundamental Theorem.} {\it Let $\mathbb{P}(V)$ and $\mathbb{P}(V')$  be projective spaces of dimensions $\geq 2$ over division rings $K$ and $K'$, respectively.}

{\it A map $\phi\colon\mathbb{P}(V)\to\mathbb{P}(V')$, such that the image is not contained in a  line, is a morphism of geometries if and only if it is induced by an injective semilinear map $\Phi\colon V\to V'$.}\medskip

Actually, the theorem was stated in more generality to characterize partial maps $\,\mathbb{P}(V)\dashrightarrow\mathbb{P}(V')\,$ induced by arbitrary semilinear maps $V\to V'$ (see Theorem \ref{5.6} below).

As a consequence, we can also characterize morphisms from a projective space $\mathbb{P}(V)$ into a vector (hence affine) space $V'$.\medskip

\noindent\textbf{Theorem.} {\it Let $\,\varphi\colon \mathbb{P}(V)\to V'\neq 0\,$ be a map whose image is not contained in an affine line of $V'$. Then $\varphi$ is a morphism of geometries if and only if it is induced by a fractional semilinear map $\,V\setminus\{0\}\to V'$.}\medskip

A map $\,\varphi\colon V\setminus\{0\}\to V'\,$ is fractional semilinear if there are semilinear maps $\,\Phi\colon V\to V'$, $\,\omega\colon V\to K'$, with the same associated morphism $K\to K'$, such that $\,\varphi(x)=\frac{1}{\omega(x)}\Phi(x)$.\medskip

Morphisms $\mathbb{A}\to\mathbb{P}\,'$, between an affine space and a projective space, are more difficult to determine. We  reduced the question to the Fundamental Theorem via the following statement, which is the main result of this paper,\medskip

\noindent\textbf{Extension Theorem.} {\it Let $\mathbb{A}$ and $\,\mathbb{P}\,'\neq *$ be an affine space and a projective space over division rings $K$ and $K'$, respectively, with $\,|K|\geq 4\,$ or $\,|K|=3=\mathrm{char}\,K'$.

Any morphism of geometries $\varphi\colon\mathbb{A}\to\mathbb{P}\,'$, whose image is not contained in a line, extends to a unique partial morphism (induced by a semilinear map) $\phi\colon\mathbb{P}\dashrightarrow\mathbb{P}\,'$, where $\mathbb{P}$ is the projective closure of $\mathbb{A}$.}
\medskip

In a forthcoming paper \cite{Sancho2}, we will use the Extension Theorem to obtain a broad generalization of the Fundamental Theorem to locally affino-projective geometries. These geometries include the so named Benz planes: The M\"{o}bius geometry of an ellipsoid, the Minkowski geometry of a ruled quadric and the Laguerre geometry of a cone.\medskip

Finally, combining the Fundamental and the Extension theorems we obtain an unpublished result of W. Zick \cite{Zick2} determining the morphisms between affine spaces,\medskip

\noindent\textbf{Theorem.} {\it Let $\mathbb{A}$ and $\mathbb{A}'\neq*$ be affine spaces over division rings $K$ and $K'$, respectively, with $|K|\geq 4$ or $|K|=3=\mathrm{char}\, K'$. Let $\varphi\colon\mathbb{A}\to\mathbb{A}'$ be a map such that the image is not contained in a line of $\mathbb{A}'$. Then, $\varphi$ is a morphism of geometries if and only if $\varphi$ is a fractional semiaffine morphism.}







\section{Geometries and morphisms}

In this section we fix the category where we shall work.\medskip

The spaces that we shall consider have different names in the literature. For the sake of brevity we shall use the name {\it geometry}, also used in \cite{FF}, \cite{Gratzer}.\medskip

\noindent\textbf{Definition.} A \textbf{geometry} is a set $X$ (whose elements are named \textbf{points}) with a family $\mathfrak{F}$ of subsets of $X$ (named  \textbf{subspaces} or \textbf{flats}), satisfying the following four axioms:\smallskip

\noindent\textbf{G1.} $\emptyset\in\mathfrak{F}$, $X\in\mathfrak{F}$ and $\{x\}\in\mathfrak{F}$ for any $x\in X$.\smallskip

The subspace $\{ x\}$ will be simply denoted by $x$.\smallskip

\noindent\textbf{G2.} Any intersection $\bigcap S_i$ of subspaces is a subspace.\smallskip

This axiom implies that the set $\mathfrak{F}$ of subspaces of $X$, with the inclusion, is a complete lattice: For any set $\{ S_i\}$ of subspaces, the infimum is
$\bigwedge S_i=\bigcap S_i$ and the supremum is $\bigvee S_i=\bigcap T_k$, where $T_k$ runs over the subspaces of $X$ containing $\bigcup S_i$.\medskip

\noindent\textbf{G3 (exchange axiom).} Given a subspace $S$ and a point $x\notin S$, there is no subspace $S'$ such that $\, S\subset S'\subset S\vee x\,$ (strict inclusions).\smallskip

That is to say, if $x,y\notin S$ then $\quad y\in S\vee x\,\Leftrightarrow\, x\in S\vee y$.\smallskip

The \textbf{closure} of a subset $A\subseteq X$ is the least subspace $\overline A$ containing $A$.\smallskip

\noindent\textbf{G4 (finitary axiom).} For any subset $A\subseteq X$ and any point $x\in \overline A$, we have $x\in \overline{\{a_1,\dots,a_r\}}$ for some finite subset  $\{a_1,\dots,a_r\}\subseteq A$.\medskip



Axioms G3 and G4 enable us to settle a beautiful dimension theory that we next resume.\medskip

\noindent\textbf{Definitions.} Given a geometry $X$, a subset $A\subseteq X$ is said to be \textbf{independent} when for any $a\in A$ we have $a\notin\overline{A\setminus a}$.

Let $A$ be a subset of a subspace $S$. If $\,\overline A=S$ we say that $A$ \textbf{generates} $S$. If $A$ is independent and generates $S$, we say that $A$ is a \textbf{basis} of $S$.





\begin{thm} Let $S$ be a subspace of a geometry $X$. Then:\smallskip

-- Any set generating $S$ contains a basis of $S$.\smallskip

-- Any independent subset of $S$ is contained in a basis of $S$.\smallskip

-- Any two bases of $S$ have equal cardinal.
\end{thm}
For a proof in a more general context, see \cite{MacLane} or \cite{FF} Chap. IV.\medskip

\noindent\textbf{Definition.} The \textbf{dimension} of a subspace is the common cardinal of its bases minus 1.

Points are just subspaces of dimension $0$. Subspaces of dimension $1$ are named \textbf{lines} and subspaces of dimension $2$, \textbf{planes}.\medskip

 Let $A=\{p_0,\dots,p_r\}\subseteq X$ be an independent finite set. Then $\overline A=p_0\vee\cdots\vee p_r$ is a subspace of dimension $r$.
In particular, any two distinct points $p_0,p_1$ lie in a unique line $\,p_0\vee p_1$, and any three non collinear points $p_0,p_1,p_2$ lie in a unique plane $\,p_0\vee p_1\vee p_2$.

\subsection*{Morphisms of geometries}

We shall come across partially defined maps. The notation $\varphi\colon X\dashrightarrow X'$ is used in order to point out that $\varphi$ is a map with domain $\text{dom}\,\varphi\subseteq X$ and image set $\text{im}\,\varphi\subseteq X'$. The \textbf{exceptional set} of $\varphi$ (with respect to $X$) is $E:=X-\text{dom}\,\varphi$. The map $\varphi$ is said to be globally defined (with respect to $X$) provided that $E$ is empty. The notation $\varphi\colon X\to X'$ is only maintained for globally defined maps.

Given any subset $S\subseteq X$, we define $\varphi(S):=\{\varphi(s):s\in S\cap\text{dom}\,\varphi\}$.  Hence $\varphi(S)=\emptyset$ exactly for $S\subseteq E$.
Given any subset $S'\subseteq X'$, we define $\varphi^{-1}(S'):=\{x\in X:\varphi(x)\in S'\text{ or }\varphi(x)=\emptyset\}=\{x\in\text{dom}\,\varphi:\varphi(x)\in S'\}\cup E$.

\begin{prop}\label{1.3} Let $\varphi\colon X\dashrightarrow X'$ be a partial map between geometries. The following statements are equivalent:\smallskip

a) The preimage of any subspace $S'\subseteq X'$ is a subspace $\varphi^{-1}(S')\subseteq X$.\smallskip

b) For any subset $A\subseteq X$ we have $\,\varphi(\,\overline A)\subseteq\overline{\varphi(A)}$, that is to say,
$$x\,\in\,\overline A\quad\Rightarrow\quad\varphi(x)\,\subseteq\,\overline{\varphi(A)}\quad .$$

c) For any finite subset $A\subseteq X$ we have $\,\varphi(\,\overline A)\subseteq\overline{\varphi(A)}$.
\end{prop}
\begin{proof} ($a\Rightarrow b$). Since $\overline{\varphi(A)}$ is a subspace of $X'$, it follows that $\varphi^{-1}\overline{\varphi(A)}$ is a subspace of $X$. Hence the inclusion $\,A\subseteq\varphi^{-1}\overline{\varphi(A)}\,$ implies $\,\overline A\subseteq\varphi^{-1}\overline{\varphi(A)}$, so that $\,\varphi(\overline A)\subseteq\overline{\varphi(A)}$. \smallskip

($b\Rightarrow a$). Let $S'\subseteq X'$ be a subspace. By ($b$) we have $\,\varphi(\overline{\varphi^{-1}S'})\subseteq\overline{\varphi(\varphi^{-1}S')}$ $\subseteq\overline{S'}=S'$, hence $\overline{\varphi^{-1}S'}\subseteq\varphi^{-1}S'$, so that $\varphi^{-1}S'$ is a subspace.\smallskip

($b\Rightarrow c$). Evident.\smallskip

($c\Rightarrow b$). Let us put $A\subseteq X$ as the union of all finite subsets: $A=\bigcup A_i$. By axiom G4, $\overline A=\bigcup\overline{A_i}$, hence $\varphi(\overline A)=\bigcup\varphi(\overline{A_i})\subseteq\bigcup\overline{\varphi(A_i)}\subseteq\overline{\varphi(A)}$.

\end{proof}

\noindent\textbf{Definition.} A partial map $\varphi\colon X\dashrightarrow X'$, between geometries, is a \textbf{partial morphism} when it fulfills the equivalent conditions of Proposition \ref{1.3}. If a partial morphism is globally defined then it is named \textbf{morphism of geometries}.

 According to Proposition \ref{1.3}c we have:

\begin{nada}\label{1.4a}  A partial map $\varphi\colon X\dashrightarrow X'$, between geometries, is a partial morphism if and only if for any finite sequence  $\,x_0,\dots,x_r\in X\,$ we have
$$x_0\,\in\, x_1\vee\dots\vee x_r\qquad\Rightarrow\qquad \varphi(x_0)\,\subseteq\,\varphi(x_1)\vee\dots\vee\varphi(x_r)\quad .$$
Recall that $\varphi(x)$ can be the empty set.
\end{nada}

The exceptional set $E=\varphi^{-1}(\emptyset)$ of a partial morphism $\varphi\colon X\dashrightarrow X'$ is a subspace of $X$.\medskip

\noindent\textbf{Definition.} We say that a geometry $X$ is \textbf{generated by lines} when the subspaces of $X$ are just the subsets $S\subseteq X$ such that
$$x_1,x_2\in S\quad \Rightarrow\quad  x_1\vee x_2\,\subseteq S\quad .$$
That is to say, in a geometry generated by lines, subspaces are just subsets containing the line joining any two different points of such subset.\medskip

\noindent\textbf{Definition.} We say that a geometry $X$ is \textbf{generated by lines and planes} when the subspaces of $X$ are just the subsets $S\subseteq X$ such that
$$x_1,x_2,x_3\in S\quad\Rightarrow\quad x_1\vee x_2\vee x_3\,\subseteq S\quad .$$
That is to say, a subset $S$ is a subspace if it contains the line joining any two different points of $S$ and the plane defined by any three non collinear points of $S$.\medskip

We shall see that projective spaces are geometries generated by lines (Proposition \ref{3.3}). Affine spaces over a division ring $K\neq\mathbb{Z}_2$ are also generated by lines (Proposition \ref{2.6}), and affine spaces over $\mathbb{Z}_2$ are generated by lines and planes.

\begin{nada}\label{1.6} Let $\varphi\colon X\dashrightarrow X'$ be a partial map between geometries, where $X$ is generated by lines. The map $\varphi$ is a partial morphism if and only if
\begin{equation}\label{eq}
x_0\,\in\, x_1\vee x_2\qquad\Rightarrow\qquad\varphi(x_0)\,\subseteq\,\varphi(x_1)\vee\varphi(x_2)
\end{equation}
for any $x_0,x_1,x_2\in X$. 

From (\ref{eq}) it easily follows (when $X$ is generated by lines) that the preimage of any subspace also is a subspace, hence $\varphi$ is a morphism. The converse holds by  \ref{1.4a}.\smallskip

Condition (\ref{eq}) is equivalent to the following three conditions: 

--  $\varphi$ transforms any three collinear points $x_0,x_1,x_2\in\text{dom}\,\varphi$ into collinear points; 

--  for any line $L\nsubseteq E\,$ the restriction $\,\varphi_{|L\!-\!E}\,$  is injective or constant;

 -- if $e\in E$ then $\,\varphi(e\vee x)=\varphi(x)\,$ for all $x\in X$.\medskip
 
 Partial maps satisfying  (\ref{eq}) are called weak linear mappings in \cite{Havlicek2}.
\end{nada}


\noindent\textbf{Definition.} A morphism of geometries $\varphi\colon X\to X'$ is an \textbf{isomorphism} if it is bijective and the inverse map $\varphi^{-1}$ also is a  morphism.\medskip

A bijective morphism need not be an isomorphism. For example, given a geometry $X$ of dimension $n\geq 2$ and a natural number $m<n$, let $X'$ be the geometry (named {\it  truncation} of $X$) with the same underlying set $X$ but such that the proper subspaces are just the subspaces of $X$ of dimension $<m$ (so that $\text{dim}\, X'=m$). The identity $X\to X'$ is a bijective morphism but not an isomorphism. There are more interesting examples; Ceccherini \cite{C} gives a bijective morphism between a 4-dimensional projective space and a non-arguesian projective plane.\medskip

 An isomorphism $\varphi\colon X\to X'$, between geometries generated by lines, is also said to be a \textbf{collineation}. That is to say, a collineation (between such geometries) is just a bijection transforming lines onto lines.

\begin{prop}\label{1.5} Let $\varphi\colon X\to X'$ be a morphism of geometries. We have:\smallskip

a) If $\varphi$ is surjective then $\,\mathrm{dim}\,X\geq\mathrm{dim}\, X'$.\smallskip

b) If $\varphi$ is surjective and $\,\mathrm{dim}\,X\leq\mathrm{dim}\, X'<\infty$ then $\varphi$ is an isomorphism.
\end{prop}
\begin{proof} $a$). If a collection $\{ x_i'=\varphi(x_i)\}$ is a basis of $X'$ (in particular, the assignment $i\mapsto x_i'$ is injective) then it is easy to see that the collection $\{ x_i\}$ is an independent set in $X$, hence $\,\mathrm{dim}\,X\geq\mathrm{dim}\, X'$.\smallskip

$b$). By $(a)$ we have $\,\mathrm{dim}\,X=\mathrm{dim}\, X'=:r$. Let $\{ x_0,\dots,x_r\}$ be a basis of $X$. Then
$$X'\,=\,\varphi(X)\,=\,\varphi(\overline{\{ x_0,\dots,x_r\}})\,\subseteq\,\overline{\{ \varphi(x_0),\dots,\varphi(x_r)\}}$$ hence $\{ \varphi(x_0),\dots,\varphi(x_r)\}$ generates $X'$ and, since $\text{dim}\, X'=r$, we conclude that $\{ \varphi(x_0),\dots,\varphi(x_r)\}$ is a basis of $X'$.

Since $\varphi$ transforms bases into bases, it also transforms independent collections into independent collections. In particular, $\varphi$ is injective, hence bijective.

We have to prove that $\varphi^{-1}$ is a morphism, that is to say, if $S$ is a subspace of $X$ then $\varphi(S)$ is a subspace of $X'$. Let $\{s_0,\dots,s_n\}$ be a  basis of $S$. We have an inclusion
$$\varphi(S)\,=\,\varphi(s_0\vee\dots\vee s_n)\,\subseteq\,\varphi(s_0)\vee\dots\vee\varphi(s_n)\,=:S'$$
This inclusion is an equality (so that $\varphi(S)=S'$ is a subspace of $X'$) because otherwise there is $x'\in S'$ such that $x:=\varphi^{-1}(x')\notin S$, so that $\{s_0,\dots,s_n,x\}$ is independent so contradicting that the image  $\{\varphi(s_0),\dots,\varphi(s_n),x'=\varphi(x)\}$ is dependent.

\end{proof}

In the somewhat different context of {\it linear spaces}, the above proposition is proved in \cite{Kreuzer}.\medskip

\noindent\textbf{Definition.} Given a geometry $X$, any subset $A\subseteq X$ is a geometry where subspaces are just subsets $S\cap A$ with $S$ a  subspace of $X$. Then we say that $A$ is a \textbf{subgeometry} of $X$.

The natural inclusion $A\hookrightarrow X$ is a morphism.\medskip

\noindent\textbf{Definition.} A morphism $\varphi\colon X\to X'$ is an \textbf{embedding} when $\varphi\colon X\to\varphi(X)$ is an isomorphism, where $\varphi(X)$ is considered as a subgeometry of $X'$.\medskip

Partial morphisms are reduced to global morphisms due to the following notion.\goodbreak\medskip

\noindent\textbf{Definition.} Let $E$ be a subspace of a geometry $X$. The \textbf{quotient space} is the quotient set
$$X/\!\!/E:=\,(X\!-\! E)/\equiv$$
where $\,\equiv\,$ is the following equivalence relation on $X-E$ :
$$x_1\equiv\,x_2\quad\Leftrightarrow\quad x_1\vee E\,=\,x_2\vee E\quad .$$

The quotient map $\,X\supseteq X\!-\! E\longrightarrow X/\!\!/E$, $x\mapsto [x]$, is a partial map $\pi\colon X\dashrightarrow X/\!\!/E$ with exceptional set $E$.

We name \textbf{subspace} of $X/\!\!/E$ any subset $\,S/\!\!/E=\pi(S)\,$ where $S$ is a subspace of $X$ containing $E$.

\begin{nada}
There exists a bijective correspondence (lattice isomorphism)
$$\begin{CD}
\left\{
\text{subspaces of } X
\text{ containing }E
\right\} @= \left\{\text{subspaces of }X/\!\!/E\right\}\\
S & \longmapsto & \pi(S)=S/\!\!/E\\
S=\pi^{-1}(S/\!\!/E) & \longleftarrow\!\shortmid & S/\!\!/E
\end{CD}$$
\end{nada}

\begin{nada}
The set $X/\!\!/E$, with the family of its subspaces, is a geometry and the quotient map $\pi\colon X\dashrightarrow X/\!\!/E$, $x\mapsto[x]$, is a partial morphism with exceptional subspace $E$. The proof is routine.
\end{nada}

\begin{lem}\label{1.7a}
Let $\varphi\colon X\dashrightarrow X'$ be a partial morphism with exceptional subspace $E$.
For all $\,x_1,x_2\in X\!-\! E\,$ we have
$$
 x_1\vee E=x_2\vee E\quad\Rightarrow\quad\varphi(x_1)=\varphi(x_2)\quad .
$$
\end{lem}
\begin{proof}  It is a consequence of the following formula
$$\overline{\varphi(x\vee E)}\,=\,\varphi(x)\vee\overline{\varphi(E)}\,=\,\varphi(x)\vee\emptyset\,=\,\varphi(x)\quad .$$

\end{proof}

The partial morphism $\pi\colon X\dashrightarrow X/\!\!/E$ has the expected universal property.
\begin{prop}\label{5.9} 
{\it Any partial morphism $\varphi\colon X\dashrightarrow X'$, with exceptional subspace $E$, uniquely factors $\varphi=\widetilde\varphi\circ\pi$ through a morphism $\,\widetilde\varphi\colon X/\!\!/E\longrightarrow X'$.
 $$\xymatrixcolsep{1pc}\xymatrixrowsep{2pc}\xymatrix{
X  \ar@{-->}[rd]_{\varphi} \ar@{-->}[rr]^{\pi} &  & X/\!\!/E \ar[ld]^{\widetilde\varphi} \\
 & X'
}$$}
\end{prop}
\begin{proof}
By Lemma \ref{1.7a} and the universal property of the quotient set with respect to an equivalence relation, the map $\varphi_{|X\!-\! E}\colon X\!-\! E\to X'$ uniquely factors $\varphi_{|X\!-\! E}=\widetilde\varphi\circ\pi_{|X\!-\! E}$ for a unique map $\,\widetilde\varphi\colon X/\!\!/E\longrightarrow X'$. We have to show that $\widetilde\varphi$ is a morphism of  geometries. Given a subspace $S'$ of $X'$, we have that $\varphi^{-1}S'$ is a subspace of $X$, hence $\pi(\varphi^{-1}S')$ is a subspace of $X/\!\!/E$. It is immediate that $\widetilde\varphi^{-1}S'=\pi(\varphi^{-1}S')$, hence $\widetilde\varphi^{-1}S'$ is a subspace of $X/\!\!/E$.

\end{proof}

\begin{nada} The domain of a partial morphism  cannot be extended in general:

Let $\varphi\colon X\dashrightarrow X'$ be a partial morphism, with exceptional subspace $E\subset X$, such that $\varphi_{|X\!-\!E}$ is not constant. If the lines of $X$ have at least three points then $\varphi$ cannot be extended to a morphism of geometries $X\supseteq U\to X'$ defined on a subgeometry $U$ strictly containing $X\!-\!E$.

The proof is an easy exercise.
\end{nada}

\section{Recollections on affine and projective spaces}

\textbf{Definition.}  An \textbf{affine space} is a set $\mathbb{A}\neq\emptyset$ (its elements are named \textbf{points}) together with a  (left) vector space $V$ over a division ring $K$ and a map  $+\colon\mathbb{A}\times V\to\mathbb{A}$, $(p,v)\mapsto p+v$, such that the following axioms are satisfied:\smallskip

(1) $\quad (p+v_1)+v_2=p+(v_1+v_2)\quad$ for all $p\in\mathbb{A}$, $v_1,v_2\in V$;\smallskip

(2) $\quad p+v=p\quad\Leftrightarrow\quad v=0\quad$ for all $p\in\mathbb{A}$, $v\in V$;\smallskip

(3) $\quad$given two points $p,\bar p\in\mathbb{A}$ there is a vector $v\in V$ (necessarily unique) such that $\bar p=p+v$.\smallskip

An affine space $(\mathbb{A},V,+)$ is simply denoted by $\mathbb{A}$.\medskip

\noindent\textbf{Definition.} A non-empty subset $S\subseteq\mathbb{A}$ is a \textbf{subspace} when it is
$$S\,=\, p+W:=\,\{p+w:w\in W\}$$
where $p\in\mathbb{A}$ is a point and $W\subseteq V$ is a vector subspace. Then $W$ is said to be the \textbf{direction} of $S$. We agree that the empty subset also is a subspace.

\begin{nada}
It is easy to check that {\it an affine space, with the family of its subspaces, is a geometry.}
\end{nada}

\noindent\textbf{Definition.} Let $V$, $V'$ be vector spaces over division rings $K$, $K'$, respectively. A map $\varphi\colon V\to V'$ is said to be \textbf{semilinear} when:\smallskip

(1) it is additive: $\varphi(v_1+v_2)=\varphi(v_1)+\varphi(v_2)\quad$ for all $v_1,v_2\in V$;\smallskip

(2) there is a ring morphism $\sigma\colon K\to K'$ such that $\,\varphi(\lambda v)=\sigma(\lambda)\varphi(v)\, $ for all $\lambda\in K$, $v\in V$.\smallskip

Remark that we do not require that the morphism $\sigma\colon K\to K'$ be surjective. If $\varphi$ is not the null map, then the morphism $\sigma$ is unique and it is said to be the \textbf{ring morphism associated} to $\varphi$.\medskip

A semilinear map $\varphi\colon V\to V'$ is said to be a \textbf{semilinear isomorphism} when $\varphi$ and $\sigma$ are bijective. In such case the inverse map  $\varphi^{-1}\colon V'\to V$ also is semilinear.\medskip

\noindent\textbf{Definition.} Let $(\mathbb{A},V,+)$ and $(\mathbb{A}',V',+)$ be affine spaces over division rings $K$ and $K'$, respectively. A map  $\varphi\colon\mathbb{A}\to\mathbb{A}'$ is a \textbf{semiaffine morphism} when there is a semilinear map $\vec\varphi\colon V\to V'$ such that
$$\varphi(p+v)\,=\, \varphi(p)+\vec\varphi(v)\qquad\qquad \forall\, p\in\mathbb{A}\,,\,v\in V\quad .$$
The semilinear map $\vec\varphi$ is unique and it is named \textbf{differential} of $\varphi$.\smallskip

A semiaffine morphism $\varphi\colon \mathbb{A}\to \mathbb{A}'$ is a \textbf{semiaffine isomorphism} or a \textbf{semiaffinity} when both $\varphi$ and the ring morphism $\sigma\colon K\to K'$ (associated to $\vec\varphi\,$) are bijective. In such case the inverse map $\varphi^{-1}\colon\mathbb{A}'\to\mathbb{A}$ also is a semiaffine isomorphism.\medskip

The prefix {\it semi} in the terms {\it semilinear, semiaffine, semiaffinity} is deleted when $K=K'$ and the associated ring morphism $\sigma\colon K\to K$ is the identity.

\begin{prop}\label{4.14} Any semiaffine map $\varphi\colon\mathbb{A}\to\mathbb{A}'$ is a morphism of geometries.
\end{prop}
\begin{proof} Let $S'\subseteq\mathbb{A}'$ a subspace. We have to show that $\varphi^{-1}(S')$ is a subspace of $\mathbb{A}$. If $\varphi^{-1}(S')=\emptyset$ it is clear. Otherwise, let $p_0\in\varphi^{-1}(S')$ and $p_0':=\varphi(p_0)\in S'$. We have $S'=p_0'+W'$ for some vector subspace $W'\subseteq V'$. It easy to check that $\varphi^{-1}(S')=p_0+{\vec\varphi\,}^{-1}(W')$, hence $\varphi^{-1}(S')$ is a subspace of $\mathbb{A}$.

\end{proof}

\begin{nada}\label{2.3} Any vector space $V$ has an underlying structure of affine space $(\mathbb{A}=V,\,V,+)$, where the map $+\colon\mathbb{A}\times V\to\mathbb{A}$ is just the addition of vectors,
$$\begin{CD}
\mathbb{A}\times V=V\times V @>{+}>> V=\mathbb{A}\quad .
\end{CD}$$

Conversely, given an affine space $(\mathbb{A},V,+)$ and a fixed point $p_0\in\mathbb{A}$ we have an affine isomorphism
$$\begin{CD}
V\quad @>{\sim}>> \quad\mathbb{A}\qquad ,\qquad v\longmapsto p_0+v\quad .
\end{CD}$$
\end{nada}


\begin{lem}\label{2.4} Let $V$ be a vector space over a division ring $K$ with $|K|\neq 2$. Let $W\subseteq V$ be a subset such that:\smallskip

a) $\,0\in W$; \smallskip

b) $\,$if $\,w_1,w_2\in W\,$ then $\,(1-t)w_1+tw_2\in W\,$  for all  $\,t\in K$ (that is to say, $\, w_1,w_2\in W\ \Rightarrow\  w_1\vee w_2\in W$).\smallskip

Then $W$ is a vector subspace.
\end{lem}
\begin{proof} It is enough to show that  $\langle w_1, w_2\rangle\subseteq W$ whenever $w_1,w_2\in W$.

Remark that if $w\in W$ then $\langle w\rangle\subseteq W$: For any $t\in K$ we have $tw=(1-t)0+tw\in W$.

Now, given $w_1,w_2\in W$, for all $x,y\in K$ we have $xw_1,yw_2\in W$, hence
$$W\,\owns\, (1-t)xw_1+tyw_2\,=\,\bar xw_1+\bar yw_2\qquad ,\quad\text{where}\quad \bar x:=(1-t)x,\, \bar y:=ty$$
Taking $t\neq 0,1$ (since $|K|\neq 2$), the values of $\bar x,\bar y$ are arbitrary, so that the vector $\bar xw_1+\bar yw_2$ is any vector of $\langle w_1, w_2\rangle$.

\end{proof}

Let $\mathbb{A}$ be an affine space over a division ring $K$. Considering it as a geometry we have:

\begin{prop}\label{2.6}
If $\,K\neq \mathbb{Z}_2$,
the affine space $\mathbb{A}$ is generated by lines.
\end{prop}
\begin{proof} Let $S\subseteq\mathbb{A}$ be a subset containing the line joining any two different points of $S$. We have to show that $S$ is a subspace. Fix $p_0\in S$ and consider the affine isomorphism $V\simeq \mathbb{A}$, $v\mapsto p_0+v$. Via this isomorphism, the subset $S$ corresponds to a subset $W\subseteq V$ fulfilling the conditions of Lemma \ref{2.4}, so that $W$ is a vector (hence affine) subspace of $V$ and, therefore, $S$ is a subspace of $\mathbb{A}$.

\end{proof}

\begin{rem}\label{2.7} Lemma \ref{2.4} and Proposition \ref{2.6} do not hold on the field $K=\mathbb{Z}_2$, but it is easy to check that in such case  $\mathbb{A}$ is generated by lines and planes.

If $|K|\geq 4$ and $E$ is  a proper subspace of $\mathbb{A}$ then the geometry $\mathbb{A}-E$ is generated by lines. We shall not use this fact.
\end{rem}

\noindent\textbf{Definition.}  The \textbf{projectivization} of a vector space $V$, or the \textbf{projective space associated} to $V$, is the set
$$\mathbb{P}(V):=\,\{\text{1-dimensional vector subspaces }\langle v\rangle\text{ of }V\}\quad .$$

Subspaces of $\mathbb{P}(V)$ are defined to be projectivizations $\mathbb{P}(W)$ of vector subspaces $W\subseteq V$. \medskip

For any pair of subspaces $\mathbb{P}(W_1),\mathbb{P}(W_2)$ of $\mathbb{P}(V)$ we have
$$\mathbb{P}(W_1)\vee\mathbb{P}(W_2)\,=\,\mathbb{P}(W_1+W_2)\qquad,\qquad \mathbb{P}(W_1)\cap\mathbb{P}(W_2)\,=\,\mathbb{P}(W_1\cap W_2)\ .$$

It is easy to check that the projective space $\mathbb{P}(V)$, with the family of its subspaces, is a geometry.

Alternatively, we may define the projective space $\mathbb{P}(V)$ as the quotient geometry $V/\!\!/\{0\}$, considering $V$ as an affine space.



\begin{prop}\label{3.3} The projective space $\mathbb{P}(V)$ is generated by lines.\end{prop}
\begin{proof} We shall use the following standard fact: Given a subspace $S\neq\emptyset$ and a point $p\notin S$ in $\mathbb{P}(V)$, the subspace $S\vee p$ is just the union of all lines joining $p$ with points of $S$. 

Let $A\subseteq\mathbb{P}$ a subset containing the line joining any two points of $A$. We have to show that $A$ is a subspace. Let $S\subseteq A$ be a maximal subspace. If $S=A$ we conclude. Otherwise there is a point $p\in A$, $p\notin S$; then
$$S\vee p\,=\, \bigcup_{s\in S}(s\vee p)\,\subseteq\, A$$
so contradicting the maximal character of $S$.

\end{proof}

\begin{rem}\label{3.4} More generally, if $K\neq\mathbb{Z}_2$ and $E$ is a proper subspace of $\mathbb{P}(V)$, then the geometry $\mathbb{P}(V)-E$ is generated by lines. We shall not use this fact.
\end{rem}

\noindent\textbf{Definition.} A semilinear map $\Phi\colon V\to V'$, $v\mapsto v'$, induces a partial map
$$\begin{CD}
\phi :\mathbb{P}(V) &\quad ---\to\quad & \mathbb{P}(V')\qquad ,\qquad \langle v\rangle\,\mapsto\,\langle v'\rangle
\end{CD}$$
with exceptional subspace $E=\mathbb{P}(\text{ker}\,\Phi)$.
We say that $\phi$ is the \textbf{partial projective morphism} associated to $\Phi$. \smallskip

Remark that $\phi$ is globally defined ($E=\emptyset$) if and only if $\Phi$ is injective ($\text{ker}\,\Phi=0$). In such case, we say that $\phi\colon \mathbb{P}(V)\to\mathbb{P}(V')$ is a \textbf{projective morphism}. A projective morphism may be not injective.

A  \textbf{projective isomorphism} is a projective morphism $\phi\colon \mathbb{P}(V)\to\mathbb{P}(V')$ associated to a semilinear isomorphism $\Phi\colon V\to V'$.\medskip

It seems to be an open problem to know if each bijective projective morphism $\mathbb{P}(V)\to\mathbb{P}(V')$, with dimensions $\geq 2$, is an isomorphism. There are examples \cite{Kreuzer2} of injective morphisms $\mathbb{P}(V)\to\mathbb{P}(V')$, preserving non-collinearity, with $\text{dim}\,\mathbb{P}(V)>\text{dim}\,\mathbb{P}(V')$.\goodbreak\medskip

If two semilinear maps $\Phi,\Phi'\colon V\to V'$ are proportional ($\Phi'=\lambda\Phi$ for some $\lambda\in K'-\{0\}$), then both induce the same partial projective  morphism $\phi=\phi'\colon\mathbb{P}(V)\dashrightarrow\mathbb{P}(V')$. Conversely,
if two semilinear maps $\Phi,\Phi'\colon V\to V'$ induce the same partial projective morphism $\phi=\phi'\colon\mathbb{P}(V)\dashrightarrow\mathbb{P}(V')$ and it is not constant, then $\Phi$ and $\Phi'$ are proportional (see \cite{F} Lemma 2.4).

\begin{prop}\label{6.2} Any partial projective morphism $\phi\colon\mathbb{P}(V)\dashrightarrow\mathbb{P}(V')$ is a partial morphism of geometries.
\end{prop}
\begin{proof} Let $\Phi\colon V\to V'$ be a corresponding semilinear map.
For any subspace $\mathbb{P}(W')\subseteq\mathbb{P}(V')$ it is easy to check that $\phi^{-1}\mathbb{P}(W')=\mathbb{P}(\Phi^{-1}W')$.

\end{proof}

Let $\mathbb{P}=\mathbb{P}(V)$ and $\mathbb{P}\,'=\mathbb{P}(V')$  be projective spaces of dimensions $\geq 2$ over division rings $K$ and $K'$, respectively. The following generalization of the classical Fundamental Theorem of Projective Geometry is due to Faure and Fr\"{o}licher \cite{FF1} and Havlicek \cite{Havlicek}; see also \cite{Havlicek2} Th. 2. A brief and elementary proof may be found in \cite{F}.

\begin{thm}[Fundamental Theorem]\label{5.6} 
 Let $\phi\colon\mathbb{P}\dashrightarrow\mathbb{P}'$ be a partial map, such that the image is not contained in a line. Then  $\phi$ is a partial projective morphism if and only if $\phi$ is a partial morphism of geometries.
 \end{thm}

\section{Morphisms $\mathbb{P}\to\mathbb{A}'$}

In this Section we will determine the morphisms of geometries $\mathbb{P}\to\mathbb{A}'$ between a projective space $\mathbb{P}$ and an affine space $\mathbb{A}'$.\medskip

\noindent\textbf{Definition.} Let $(\mathbb{A},V,+)$ be an affine space over a division ring $K$. We name \textbf{vectorial extension} of $\mathbb{A}$ to a (left) $K$-vector space $\mathbb{E}$ together with an injective affine morphism $j\colon\mathbb{A}\hookrightarrow\mathbb{E}$ defining an affinity between $\mathbb{A}$ and an affine  hyperplane of $\mathbb{E}$ not containing $0$.\smallskip

The differential of $j\colon\mathbb{A}\hookrightarrow\mathbb{E}$ is an injective $K$-linear map $\vec j\colon V\hookrightarrow\mathbb{E}$.

\begin{nada}\label{2.8}{\it Comments.} 
a). It is easy to construct a vectorial extension: Fix a point $p_0\in\mathbb{A}$, put formally $\mathbb{E}:=Kp_0\oplus V$ and let us  consider the natural inclusion
$$\begin{CD}
\mathbb{A}=p_0+V @>{j}>> Kp_0\oplus V=\mathbb{E}\quad ,\qquad p_0+v\longmapsto 1\cdot p_0+v
\end{CD}$$

Its differential $\vec j$ is the natural inclusion
$$\begin{CD}
V @>{\vec j}>> Kp_0\oplus V=\mathbb{E}\quad ,\qquad v\longmapsto 0\cdot p_0+v
\end{CD}$$

In conclusion, we may identify $\mathbb{A}$ with an affine hyperplane of $\mathbb{E}$ and $V$ with a vectorial hyperplane of $\mathbb{E}$, both hyperplanes being parallel.\smallskip

b). The vectorial extension is essentially unique: Given any two vectorial extensions $j\colon\mathbb{A}\to\mathbb{E}$ and $j'\colon\mathbb{A}\to\mathbb{E}'$, there is a unique $K$-linear isomorphism $\Phi\colon\mathbb{E}\to\mathbb{E}'$ such that $j'=\Phi\circ j$.\smallskip

c). {\it Universal property of vectorial extensions:}
Any semiaffine morphism $\phi\colon\mathbb{A}\to\mathbb{A}'$ uniquely extends to a semilinear map $\Phi\colon\mathbb{E}\to\mathbb{E}'$ between their vectorial extensions (with the same associated ring morphism $K\to K'$).

The proof is easy:
Upon fixing points $p_0\in\mathbb{A}$, $p_0'=\phi(p_0)\in\mathbb{A}'$, let us consider the vectorial extensions $\mathbb{A}=p_0+V\hookrightarrow Kp_0\oplus V=\mathbb{E}$, $\,\mathbb{A}'=p_0'+V'\hookrightarrow K'p_0'\oplus V'=\mathbb{E}'$. Then the semilinear map $\,\Phi\colon Kp_0\oplus V\to K'p_0'\oplus V'\,$ is defined by the formulae $\Phi(p_0)=p_0'$, $\Phi_{|V}=\vec\phi$.
\end{nada}

 Let $(\mathbb{A},V,+)$ be an affine space over a division ring $K$. Let us recall the notations of \ref{2.8}a,
 $$\begin{CD}
\mathbb{A}\,=\,p_0+V @>{j}>> Kp_0\oplus V\,=\,\mathbb{E}\quad ,&\qquad p_0+v&\quad\longmapsto\quad& 1\cdot p_0+v\\
 V @>{\vec j}>> Kp_0\oplus V\,=\,\mathbb{E}\quad ,&\qquad v&\quad\longmapsto\quad&  0\cdot p_0+v
 \end{CD}$$
 
 \smallskip

 \noindent\textbf{Definition.} The projective space $\mathbb{P}=\mathbb{P}(\mathbb{E})$ is said to be the \textbf{projective closure} of $\mathbb{A}$. The hyperplane  $\,H:=\mathbb{P}(V)\,$ of $\mathbb{P}(\mathbb{E})$ is said to be the \textbf{hyperplane at infinity}.

 \begin{thm}\label{4.1} The composition
 $$\begin{CD}
 \mathbb{A} @>{j}>> \mathbb{E}\!-\{0\} @>{\pi}>> \mathbb{P}(\mathbb{E})=\mathbb{P}\\
 p_0+v & \longmapsto & p_0+v &\longmapsto & \langle p_0+v\rangle
 \end{CD}$$
defines a canonical isomorphism of geometries $\,\mathbb{A} =\!=\mathbb{P}\!-\! H$.
 \end{thm}

See \cite{Berger}, Th. 5.1.3 and Prop. 5.3.2, for a proof (with a different language). 

\begin{nada} {\it Comment.} Recall that any semiaffine morphism $\phi\colon\mathbb{A}\to\mathbb{A}'$ uniquely extends to a semilinear map $\Phi\colon\mathbb{E}\to\mathbb{E}'$ between their vectorial extensions (see \ref{2.8}c).
Now, this semilinear map defines a partial projective morphism $\phi\colon\mathbb{P}=\mathbb{P}(\mathbb{E})\dashrightarrow\mathbb{P}(\mathbb{E}')=\mathbb{P}\,'$. Therefore, 

 {\it Any semiaffine morphism $\phi\colon\mathbb{A}\to\mathbb{A}\,'$ uniquely extends to a partial projective morphism between their projective closures.}
 $$\xymatrixcolsep{3pc}\xymatrixrowsep{0.5pc}\xymatrix{
 \mathbb{A} \ar[r]^\phi & \mathbb{A}'\\
 \bigcap & \bigcap \\
 \mathbb{P} \ar@{-->}[r]^{\phi} & \mathbb{P}\,'
 }$$
\end{nada}

\subsection*{Fractional semilinear maps}

Let $V$, $V'$ be vector spaces over division rings $K$, $K'$, respectively. We also consider them as affine spaces, hence geometries.\medskip

\noindent\textbf{Definition.}   A map $\varphi\colon V-\{0\}\longrightarrow V'$ is said to be a \textbf{fractional semilinear map} if there exist  semilinear maps $\Psi\colon V\to V'$ and $\omega\colon V\to K'$, with the same associated ring morphism $ K\to K'$, such that 
$$\text{ker}\,\omega\,=\,0\quad\text{ and }\quad \varphi(v)\,=\,\frac{1}{\omega(v)}\Psi(v)\qquad \forall v\in V-\{0\}\quad .$$
\smallskip

Note that $\varphi(\lambda v)=\varphi(v)$ for all $v\in V\!-\!\{0\}$, $\lambda\in K\!-\!\{0\}$, hence  $\varphi$ induces a map $\,\bar\varphi\colon\mathbb{P}(V)\to V'$, $\bar\varphi\langle v\rangle:=\varphi(v)$.

\begin{prop} Let $\varphi\colon V-\{0\}\longrightarrow V'$ be a fractional semilinear map. The induced map $\,\bar\varphi\colon\mathbb{P}(V)\to V'$ is a morphism of geometries.
\end{prop}
\begin{proof} Since $\mathbb{P}(V)$ is generated by lines, we have to prove the property
$$\langle v_0\rangle\in\,\langle v_1\rangle\vee\langle v_2\rangle \quad\Rightarrow\quad \varphi(v_0)\in\,\varphi(v_1)\vee\varphi(v_2)\quad .$$

If $\langle v_0\rangle\in\,\langle v_1\rangle\vee\langle v_2\rangle$ then $v_0\in\langle v_1,v_2\rangle$ hence $v_0=\lambda_1v_1+\lambda_2v_2$ for some $\lambda_1,\lambda_2\in K$. Now, a direct computation gives
$$\varphi(v_0)\,=\,\varphi(\lambda_1v_1+\lambda_2v_2)\,=\,\frac{\omega(\lambda_1v_1)}{\omega(\lambda_1v_1+\lambda_2v_2)}\varphi(v_1)+\frac{\omega(\lambda_2v_2)}{\omega(\lambda_1v_1+\lambda_2v_2)}\varphi(v_2)$$
where we use the convention $\frac{a}{b}:=b^{-1}a$. Denote $\mu_i=\frac{\omega(\lambda_iv_i)}{\omega(\lambda_1v_1+\lambda_2v_2)}\in K'$, $i=1,2$. The formulae
$$\varphi(v_0)\,=\,\mu_1\varphi(v_1)+\mu_2\varphi(v_2)\qquad ,\qquad \mu_1+\mu_2=1$$
imply that $\,\varphi(v_0)\in\,\varphi(v_1)\vee\varphi(v_2)$.

\end{proof}

\begin{thm}\label{4.4} Let $\,\bar\varphi\colon \mathbb{P}(V)\to V'\neq 0\,$ be a map whose image is not contained in an affine line of $V'$. Then $\bar\varphi$ is a morphism of geometries if and only if it is induced by a fractional semilinear map $\,\varphi\colon V-\{0\}\to V'$.
\end{thm}
\begin{proof}  ($\Rightarrow$). The vectorial extension of $V'$, as an affine space, is
$$\begin{CD}
 V' @= 1+V' & \,\subset\, K'\oplus V'\,=:\mathbb{E}'
\\ v' &\longmapsto & 1+v'
\end{CD}$$

The projective embedding of $V'$, as an affine space, is
$$\begin{CD}
V' @= 1+V' & \quad\subset\quad & \mathbb{P}(\mathbb{E}')\\
v' &\longmapsto & 1+v' &\longmapsto & \langle 1+v'\rangle
\end{CD}$$ 

Let $\varphi\colon V-\{0\}\to V'$, $\varphi(v):=\bar\varphi\langle v\rangle$. By Fundamental Theorem \ref{5.6}, the morphism of geometries 
$$\begin{CD}
\mathbb{P}(V) @>{\bar\varphi}>> V' &\quad \subset\quad &\mathbb{P}(\mathbb{E}')\\
\langle v\rangle &\longmapsto & \bar\varphi\langle v\rangle=\varphi(v) & \longmapsto & \langle\,1+\varphi(v)\,\rangle
\end{CD}$$ is induced by a semilinear map 
$$\begin{CD}
V @>{(\omega,\Psi)}>>  K'\oplus V' & \,=\,\mathbb{E}'\\
v &\longmapsto & \omega(v)+\Psi(v)
\end{CD}$$ that is to say,
$$\begin{CD}
\langle\,1+\varphi(v)\,\rangle @= \langle\,\omega(v)+\Psi(v)\,\rangle\,=\, \langle\, 1+\frac{1}{\omega(v)}\Psi(v)\,\rangle
\end{CD}$$
hence $\,\varphi(v)=\frac{1}{\omega(v)}\Psi(v)\,$ is a fractional semilinear map.\smallskip

($\Leftarrow$). It is the previous proposition.
\end{proof}

\begin{rem} There are no fractional semilinear maps $\varphi=\frac{1}{\omega}\Psi\colon V-\{0\}\to V'$ when $V,V'$ are real vector spaces with $\text{dim}\,V\geq 2$. 
\end{rem}
\begin{proof} Since any ring morphism $\mathbb{R}\to\mathbb{R}$ is the identity, the map $\,\omega\colon V\to\mathbb{R}\,$ is $\mathbb{R}$-linear. Then the condition $\,\text{ker}\,\omega= 0\,$ is impossible when $\text{dim}\, V\geq 2$.

\end{proof}

\subsection*{Partial morphisms $\mathbb{A}\dashrightarrow\mathbb{A}'$}

Theorem \ref{4.4} can be used to determine the (non global) partial morphisms of geometries $\mathbb{A}\dashrightarrow\mathbb{A}'$ between affine spaces. The result is stated more clearly in terms of vector spaces.

\begin{thm} Let $V$ and $V'\neq 0$ vector spaces over division rings $K$ and $K'$, respectively. Let $W\subset V$ be a vector subspace and let $\varphi\colon V\!-\!W\longrightarrow V'$ be a map whose image is not contained in an affine line of $V'$. The following conditions are equivalent:\smallskip

a) $\,\varphi\colon V\dashrightarrow V'$ is a partial morphism of geometries with exceptional subspace $W$.\smallskip

b) There exist semilinear maps $\Psi\colon V\to V'$, $\omega\colon V\to K'$, with the same associated ring morphism $K\to K'$, such that
$$W\,=\,\mathrm{ker}\,\omega\,\subseteq\,\mathrm{ker}\,\Psi\quad \text{ and } \quad \varphi(v)\,=\,\frac{1}{\omega(v)}\Psi(v)\quad \forall\, v\in V\!-\! W\quad .$$
\end{thm}
\begin{proof}  By Proposition \ref{5.9}, condition $(a)$ means  that $\varphi$ factors through a morphism $\bar\varphi\colon V/\!\!/W\to V'$. Now, it is routine to check that $V/\!\!/W=\mathbb{P}(V/W)$. We conclude applying Theorem\ref{4.4} to $\,\bar\varphi\colon V/\!\!/W=\mathbb{P}(V/W)\longrightarrow V'$.

\end{proof}

\section{Global morphisms $\mathbb{A}\to\mathbb{P}\,'$}


\subsection*{Fractional semiaffine morphisms}

Let $(\mathbb{A},V,+)$ and $(\mathbb{A}',V',+)$ be affine spaces over division rings $K$ and $K'$, respectively. Let $\mathbb{A}'\hookrightarrow\mathbb{E}'$ be the vectorial extension of $\mathbb{A}'$.\medskip

\noindent\textbf{Definition.}   A map $\varphi\colon\mathbb{A}\to \mathbb{A}'\subset\mathbb{E}'$ is said to be a \textbf{fractional semiaffine morphism} if there exist  semiaffine morphisms $\Psi\colon\mathbb{A}\to\mathbb{E}'$ and $\omega\colon\mathbb{A}\to K'$, with the same associated ring morphism $ K\to K'$, such that 
$$\omega(x)\,\neq\,0\quad\text{ and }\quad  \varphi(x)\,=\,\frac{1}{\omega(x)}\Psi(x)\qquad \forall x\in\mathbb{A}\quad .$$

\smallskip

\noindent{\it Comments}

-- Let $\mathbb{A}'\neq *$. Let us assume that either $(a)$ $\text{char}\, K\neq 2$ and $\varphi$ is not constant  or $(b)$ $\text{char}\, K=2$ and the image of $\varphi$ is not contained in a line of $\mathbb{A}'$. In both cases the pair $\Psi,\omega$ is unique except for  a common scalar factor $\lambda'\in K'$. We shall not use this fact.\smallskip

-- Let us describe the equations of a fractional semiaffine morphism $\varphi=\frac{1}{\omega}\Psi\colon\mathbb{A}_n\to\mathbb{A}_m'$, between affine spaces of finite dimensions. If $(x_1,\dots,x_n)$ and $(y_1,\dots,y_m)$ are affine coordinates in $\mathbb{A}_n$ and $\mathbb{A}'_m$, respectively, then the equations of $\varphi$ have the form
$$\left\{\begin{aligned}
y_1\,&=\,\frac{1}{c_0+x_1'c_1+\cdots +x_n'c_n}(a_{10}+x_1'a_{11}+\cdots+x_n'a_{1n})\\
\vdots\, &   \\
y_m\,&=\,\frac{1}{c_0+x_1'c_1+\cdots +x_n'c_n}(a_{m0}+x_1'a_{m1}+\cdots+x_n'a_{mn})
\end{aligned}\right.$$
where $a_{ij},c_k\in K'$ and $K\to K'$, $x\mapsto x'$, is the common homomorphism associated to $\Psi$ and $\omega$.\smallskip

--  A map $\varphi\colon V\to V'$, between vector spaces, is a fractional semiaffine morphism if and only if  there exist  semiaffine morphisms $\Phi\colon V\to V'$ and $\omega\colon V\to K'$, with the same associated ring morphism $ K\to K'$, such that 
$$\omega(v)\,\neq\,0\quad\text{ and }\quad  \varphi(v)\,=\,\frac{1}{\omega(v)}\Phi(v)\qquad \forall v\in V$$

\begin{rem}\label{5.1a} {\it Any fractional semiaffine morphism $\varphi=\frac{1}{\omega}\Psi\colon\mathbb{A}\to\mathbb{A}'$, between real affine spaces, is an affine map.} 
\end{rem}
\begin{proof} Since any homomorphism $\mathbb{R}\to\mathbb{R}$ is the identity, the maps $\Psi,\omega$ are affine morphisms. The condition $\,\omega(x)\neq 0$ for any $x\in\mathbb{A}\,$ implies that the affine function $\omega\colon\mathbb{A}\to\mathbb{R}$ is constant, hence $\varphi=\frac{1}{\omega}\Psi$ is an affine map.

\end{proof}

Let $\mathbb{A}\subset\mathbb{E}$ be the vectorial extension of an affine space. A linear combination $\sum_1^r\lambda_ip_i$ of a finite collection of points of $\mathbb{A}$ has sense as a vector in $\mathbb{E}$. An standard fact says that the affine subspace $p_1\vee\cdots\vee p_r\subseteq\mathbb{A}$ coincides with the following set of linear combinations (named \textbf{affine combinations}):
$$p_1\vee\cdots\vee p_r\,=\,\left\{\sum\lambda_ip_i:\sum\lambda_i=1\right\}\quad .$$

\begin{prop}\label{4.15}
Any fractional semiaffine morphism $\varphi\colon\mathbb{A}\to\mathbb{A}'$, $\varphi(x)=\frac{1}{\omega(x)}\Psi(x)$, is a morphism of geometries.
\end{prop}
\begin{proof} Let $K\to K'$, $\lambda\mapsto\lambda'$, be the common associated ring morphism of $\Psi$ and $\omega$. We have to prove
$$p_0\,\in\, p_1\vee \cdots\vee p_r\qquad\Rightarrow\qquad\varphi(p_0)\,\in\,\varphi(p_1)\vee\cdots\vee\varphi(p_r)$$
for any $p_0,p_1,\dots,p_r\in\mathbb{A}$.

If  $p_0\in p_1\vee \cdots\vee p_r$ then $p_0=\sum_1^r\lambda_ip_i$ with $\sum_1^r\lambda_i=1$, hence
$$\varphi(p_0)\,=\,\varphi\left(\sum\lambda_ip_i\right)\,=\,\frac{1}{\omega\left(\sum\lambda_ip_i\right)}\Psi\left(\sum\lambda_ip_i\right)\,=\,\sum\frac{\lambda_i'}{\omega\left(\sum\lambda_ip_i\right)}\Psi(p_i)$$
$$=\,\sum\frac{\lambda_i'\,\omega(p_i)}{\omega\left(\sum\lambda_ip_i\right)}\varphi(p_i)
$$
where the last term is an affine combination:
$$\sum\frac{\lambda_i'\,\omega(p_i)}{\omega\left(\sum\lambda_ip_i\right)}\,=\,\frac{\sum \lambda_i'\omega(p_i)}{\omega\left(\sum\lambda_ip_i\right)}\,=\, \frac{\omega\left(\sum\lambda_ip_i\right)}{\omega\left(\sum\lambda_ip_i\right)}\,=\,1$$
Therefore, $\varphi(p_0)\,\in\,\varphi(p_1)\vee\cdots\vee\varphi(p_r)$.

\end{proof}

\subsection*{Projective extension of morphisms $\mathbb{A}\to\mathbb{P}\,'$}

We are going to determine the morphisms $\mathbb{A}\to\mathbb{P}\,'$, between an affine space and a projective space. We will reduce the question to the Fundamental Theorem, proving that each morphism $\mathbb{A}\to\mathbb{P}\,'$ extends to a partial morphism $\mathbb{P}\dashrightarrow\mathbb{P} \,'$, where $\mathbb{P}$ is the projective closure of $\mathbb{A}$.

From now on, we consider an affine space $(\mathbb{A},V)$ over a division ring $K$, its projective closure $\mathbb{P}=\mathbb{P}(\mathbb{E})$ and the hyperplane at infinity $H=\mathbb{P}(V)\subset\mathbb{P}$. According to Theorem \ref{4.1} we have an identification  $\mathbb{A}=\mathbb{P}\!-\! H$ defined by the composition
$$\begin{CD}
 \mathbb{A} @>{j}>> \mathbb{E}\!-\{0\} @>{\pi}>> \mathbb{P}(\mathbb{E})=\mathbb{P}\\
 p_0+v & \longmapsto & p_0+v &\longmapsto & \langle p_0+v\rangle
 \end{CD}$$\smallskip

We frequently use that the closure in $\mathbb{P}$ of an affine line $L\subseteq\mathbb{A}$, with direction $\langle v\rangle\subseteq V$, is a projective line meeting the hyperplane at infinity $H$ at a point $\,p=\langle v\rangle$, named \textbf{point at infinity} of the affine line $L$. Therefore,
{\it two affine lines are parallel if and only if both have the same point at infinity.}\medskip

Let us consider another projective space $\mathbb{P}\,'\neq *$ over a division ring $K'$.

We assume that $|K|\geq 4$ or $|K|=3=\text{char}\,K'$.\medskip

Finally, let us fix a morphism of geometries
$$\begin{CD}
\varphi: \mathbb{A} @>>> \mathbb{P}\,'\end{CD}$$ such that the image is not contained in a line.\medskip

\noindent\textbf{Notation.} Since $\varphi$ is a morphism, for any line $L$ of $\mathbb{A}$ we have that $\varphi_{|L}$ is constant or $\varphi_{|L}$ is injective and  $\varphi(L)$ is contained in a unique line of $\mathbb{P}\,'$, denoted by $L'$.

\begin{lem}\label{4.3} Given a point $p\in H$, let $\{ L_i\}$ be the family of parallel lines in $\mathbb{A}$ whose point at infinity is $p$. Let us consider the family $\{L_j'\}$ of lines in $\mathbb{P}\,'$, where $L_j$ runs over the lines in $\{L_i\}$ such that $\varphi_{|L_j}$ is injective. One of the following alternatives holds:\smallskip

i) $\varphi$ is constant on any line in the family $\{L_i\}$ (that is to say, $\{ L_j'\}=\emptyset$).\smallskip

ii) The union of the lines in the family $\{L_j'\}$ contains $\varphi(\mathbb{A})$.
\end{lem}
\begin{proof} If $(i)$ does not hold, then there is a line $L_0$ in the family $\{L_i\}$ such that $\varphi_{|L_0}$ is injective. Let us prove $(ii)$. Given $x\in \mathbb{A}$, let us see that $\varphi(x)\in L_j'$ for some index $j$. If $\varphi(x)\in L_0'$ it is obvious. Otherwise, $\varphi(x)$ and $\varphi(L_0)$ are not collinear. Consider the plane $P=L_0\vee x$ in the geometry $\mathbb{A}$. Since $\text{dim}\,\varphi(P)=2$, we have that $\varphi\colon P\to\varphi(P)$ is an isomorphism by Proposition \ref{1.5}b. Now let $L_j$ be the line in the plane $P$ parallel to $L_0$ passing through $x$. Since $\varphi_{|P}$ is injective, so is $\varphi_{|L_j}$. Finally, since $x\in L_j$ we have $\varphi(x)\in\varphi(L_j)\subseteq L_j'$.

\end{proof}

\begin{nada}\label{4.4b} Points $p\in H$ where $(i)$ holds are said to be \textbf{exceptional}.\end{nada}

\begin{lem}\label{4.4a} The set $E\subseteq H$ of exceptional points is a subspace.
\end{lem}
\begin{proof} Let $p_1,p_2\in H$ be two different exceptional points. We must show that any point $p\in p_1\vee p_2\subseteq H$ also is exceptional. Let $L$ be an affine line in $\mathbb{A}$ with the point $p$ at infinity and let $x\in L$. Let us consider the lines $L_1,L_2$ of $\mathbb{A}$ passing through $x$ with the points $p_1,p_2$ at infinity respectively. Then $L_1\vee L_2$ is a plane of $\mathbb{A}$ containing the line $L$. Since $p_1,p_2$ are exceptional, we have $\varphi(L_1)=\varphi(x)$, $\varphi(L_2)=\varphi(x)$. Hence $\varphi(L)\subseteq\varphi(L_1\vee L_2)\subseteq\varphi(L_1)\vee\varphi(L_2)=\varphi(x)\vee\varphi(x)=\varphi(x)$, so that $\varphi_{|L}$ is constant.

\end{proof}

\begin{lem}\label{7.6} If $p\in H$ is not an exceptional point (according to \ref{4.4b}), then the lines of $\mathbb{P}\,'$ in the family $\{L_j'\}$ meet at a unique point $p'$ of $\,\mathbb{P}\,'$.
\end{lem}
\begin{proof} Let $L_1',L_2'$ be two different lines in the family $\{L_j'\}$. Since the lines $L_1,L_2$ of $\mathbb{A}$ are coplanar, the lines $L_1',L_2'$ of $\mathbb{P}\,'$ so are coplanar, hence they intersect each other. Now we must prove that any three different lines $L_1',L_2',L_3'$ meet at a point.\smallskip

Case 1: The parallel lines $L_1,L_2,L_3$ lie in a plane $\mathbb{A}_2$. In this case $L_1',L_2',L_3'$ are coplanar and $\text{dim}\,\varphi(\mathbb{A}_2)=2$, hence $\varphi\colon \mathbb{A}_2\to\varphi(\mathbb{A}_2)$ is an isomorphism, that is to say, $\varphi$ induces an embedding of $\mathbb{A}_2$ into a projective plane  $\mathbb{P}_2'\subseteq\mathbb{P}\,'$. \smallskip

1a. If $|K|\geq 4$, there is a forth line $L\subset\mathbb{A}_2$, different and parallel to the other three lines (figure 1 below). Let us pick different points $x_1,\bar x_1\in L_1$, $z_{12},z_{13}\in L$ and let us consider the points $x_2=L_2\cap (x_1\vee z_{12})$, $\bar x_2=L_2\cap(\bar x_1\vee z_{12})$, $x_3=L_3\cap (x_1\vee z_{13})$, 
$\bar x_3=L_3\cap(\bar x_1\vee z_{13})$. The vertices of the triangles $x_1x_2x_3$ and $\bar x_1\bar x_2\bar x_3$ lie in the parallel lines $L_1,L_2,L_3$, hence by  Desargues's theorem the point $z_{23}:=(x_2\vee x_3)\cap(\bar x_2\vee \bar x_3)$ is collinear with $z_{12}$ and $z_{13}$, that is to say, $z_{23}\in L$.
\begin{figure}[h]
  \centering
  \includegraphics[height=50mm]{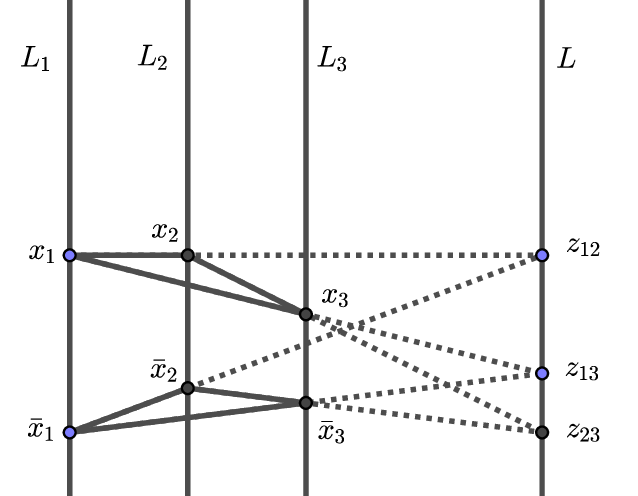}
  \caption{}
\end{figure}

Applying to these triangles the embedding $\varphi\colon\mathbb{A}_2\to\mathbb{P}_2'$ we obtain two triangles in $\mathbb{P}_2'$, with vertices in the lines 
$L_1',L_2',L_3'$, such that the corresponding sides meet at points in the line $L'$. By Desargues's theorem, the lines $L_1',L_2',L_3'$ are concurrent. We take this  argument from \cite{Rigby} Th 6.2.\medskip

1b. If $|K|=3=\text{char}\,K'$ then $\mathbb{A}_2$ is an affine plane over $K=\mathbb{Z}_3$, with nine points and 12 lines. Therefore, the image of the embedding 
$\varphi\colon\mathbb{A}_2\to\mathbb{P}_2'$ has 9 points as follows (figure 2).
\begin{figure}[h]
  \centering
  \includegraphics[height=60mm]{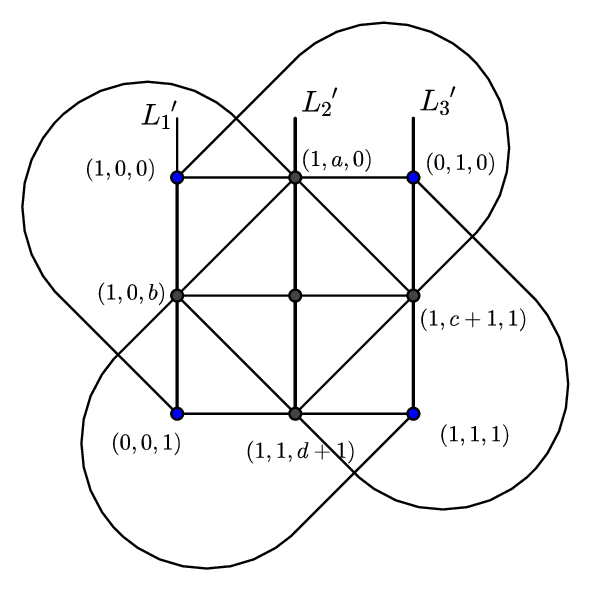}
  \caption{}
\end{figure}

Let us fix the projective reference in $\mathbb{P}_2'$ defined by the four vertices of the big square. Then the coordinates of the vertices of the inscribed rhombus are those in the figure, with $a,b,c,d\in K'-\{0\}$. Since the two vertices of each side of the rhombus are collinear with the opposite vertex of the square, a direct calculation (using that $\text{char}\, K'=3$) shows that:
 $$a=b=2\qquad ,\qquad c=d=1\quad .$$
Now we check, using that $\text{char}\, K'=3$, that the lines $L_1',L_2',L_3'$ are concurrent at the point $(1,0,1)$. We take this argument from \cite{OS} Th.2.\goodbreak\smallskip

Case 2: The parallel lines $L_1,L_2,L_3$ generate a $3$-space $\mathbb{A}_3$. We argue in a similar way to case 1a; the triangles being constructed as follows. Let $P\subset\mathbb{A}_3$ 
be a plane not parallel to the lines $L_1,L_2,L_3$, so that it meet them at the vertices of a triangle. Let $L\subset P$ be a line not parallel to any side of the former  triangle (it exists because $|K|>2$). Take another plane $\bar P\subset\mathbb{A}_3$, with the same property than $P$, intersecting $P$ at $L$. This plane $\bar P$ meets the lines $L_1,L_2,L_3$ at the vertices of a second triangle. By construction, the vertices of both triangles lie in the lines $L_1,L_2,L_3$ and the corresponding sides meet at points of $L$. 
We conclude as in case 1a.

\end{proof}

The above lemma let us extend the morphism $\varphi\colon\mathbb{A}\to\mathbb{P}\,'$ to a partial map $\phi\colon\mathbb{P}\dashrightarrow\mathbb{P}\,'$ taking $\phi(p):=p'$ for any non exceptional point $p\in H$.

\begin{lem}\label{8.4} Let $p\in H$ be a non exceptional point and let $L_0$ be a line of $\mathbb{A}$ with the point $p$ at infinity.\smallskip

a) If $\varphi_{|L_0}$ is constant then $\,\varphi(L_0)=p'$.\smallskip

b) If $\varphi_{|L_0}$ is injective then $\,p'\notin\varphi(L_0)$.
\end{lem}
\begin{proof} $a)$. Put $x':=\varphi(L_0)$. Given a line $L_j'$ in the family $\{L_j'\}$, let us consider the plane $P$ of $\mathbb{A}$ containing the parallel lines $L_0$ and $L_j$. Since $\varphi(L_0)=x'$ we have $\text{dim}\,\varphi(P)\leq 1$, so that $x'=\varphi(L_0)$ and $\varphi(L_j)$ are collinear, that is to say, $x'\in L_j'$. 
Therefore, $x'$ is the point where the lines in the family $\{L_j'\}$ intersect, that is to say, $x'=p'$.\smallskip

$b)$. If $p'\in\varphi(L_0)$ we obtain a contradiction. Let $\varphi(x)\in\varphi(\mathbb{A})$ be a point, not collinear with $\varphi(L_0)$. Let us consider the plane $P\subseteq\mathbb{A}$ containing $L_0$ and $x$; let $L_j\subset P$ be the parallel line to $L_0$ through $x$. Remark that $\text{dim}\, \varphi(P)=2$ because  $\varphi(P)\supset \varphi(L_0),\varphi(x)$, so that $\varphi\colon P\to\varphi(P)$ is an isomorphism. It follows that $\varphi(L_j)$ is a line of $\varphi(P)$, hence $\varphi(L_j)=L_j'\cap\varphi(P)\owns p'$. Therefore, $\varphi(L_0)$ and $\varphi(L_j)$ meet at $p'$, so contradicting that $\varphi\colon P\to\varphi(P)$ is an isomorphism and $L_0\cap L_j=\emptyset$.

\end{proof}

\begin{lem}\label{8.5} For any $\,p_0,p_1,p_2\in H-E\,$ we have
$$p_0\,\in\,p_1\vee p_2\quad\Rightarrow\quad p_0'\,\in\,p_1'\vee p_2'\quad .$$
\end{lem}
\begin{proof} 
Case $p_1'=p_2'=:\bar p$. We have to show that $p_0'=\bar p$, that is to say, for any line $L_0\subset\mathbb{A}$ with the point $p_0$ at infinity and  such that $\varphi_{|L_0}$ is injective, we have $\bar p\in L_0'$. Given $x_0\in L_0$, let $L_1,L_2\subset\mathbb{A}$ be the lines through $x_0$ with the points at 
infinity $p_1,p_2$, respectively. If $\varphi_{|L_1}$ is constant, then $\bar p=p_1'=\varphi(L_1)\owns\varphi(x_0)$, so that $\bar p=\varphi(x_0)\in\varphi(L_0)\subseteq L_0'$. It is the same if $\varphi_{|L_2}$ is constant. Finally, if $\varphi$ is injective on $L_1$ and $L_2$, then $L_1'=L_2'$ because both lines pass through $\varphi(x_0)$ and $\bar p$; the inclusion $L_0\subset L_1\vee L_2$ implies  $\varphi(L_0)\subseteq\varphi(L_1)\vee\varphi(L_2)\subseteq (L_1'=L_2')$, hence $L_0'=L_1'=L_2'\owns \bar p$.\smallskip

Case $p_1'\neq p_2'$. We may assume that $p_0\neq p_1,p_2$. Let us consider the plane $P\subseteq\mathbb{A}$ with the line $p_1\vee p_2$ at infinity and passing through a point $x_0\in\mathbb{A}$ such that $\varphi(x_0)\notin p_1'\vee p_2'$.

Let us see that $\text{dim}\,\varphi(P)=2$ (so that $\varphi\colon P\to\mathbb{P}\,'$ is an embedding).
Let $L_i\subset P$ be the line through $x_0$ with the point $p_i$ at infinity ($i=1,2$). Remark that $\varphi_{|L_i}$ is injective since otherwise  $p_i'=\varphi(L_i)\owns \varphi(x_0)$, against the choice of $x_0$.
If $\text{dim}\,\varphi(P)\leq 1$ then the subsets $\varphi(L_1),\varphi(L_2)\subseteq\varphi(P)$ lie in a line, hence $L_1'=L_2'$, so that 
$\varphi(x_0),p_1',p_2'$ are collinear, against the choice of $x_0$.

Now let us consider in the affine plane $P$ two perspective triangles such that the corresponding sides are parallel with the points $p_0,p_1,p_2$ at infinity 
(such triangles exist because $P$ has enough points when $|K|>2$). Applying the embedding $\varphi\colon P\to\mathbb{P}\,'$ to such configuration, we obtain a pair 
of perspective triangles in $\mathbb{P}\,'$ such that the corresponding sides intersect at the points $p_0',p_1',p_2'$, hence these points are collinear by Desargues's theorem.

 \end{proof}
 
\begin{lem}\label{6.6} Let $p_1,p_2\in H-E$ be two different points. Then
$$(p_1\vee p_2)\cap E\neq\emptyset\qquad\Rightarrow\qquad p_1'=p_2'\quad .$$
\end{lem}
\begin{proof} Put $p=(p_1\vee p_2)\cap E$. Let  $\{I_i\}$ be the family of parallel lines in $\mathbb{A}$ with point at infinity $p_1$ and let $\{I_j'\}$ be the corresponding family of lines in $\mathbb{P}\,'$, concurrent at $p_1'$. Given one of such lines $I'$, let us consider in $\mathbb{A}$ a plane $P\supset I$, with the  line at infinity $p_1\vee p_2$. Let $J,L\subset P$ be lines with the points at infinity $p_2,p$ respectively. Since $p$ is exceptional we have that $\varphi_{|L}$ is constant, hence $\text{dim}\,\varphi(P)\leq 1$ and therefore $\varphi(P)\subseteq I'$. We have $\varphi(J)\subseteq\varphi(P)\subseteq I'$; hence, if  $\varphi_{|J}$ is injective, we have $J'=I'$ so that $p_2'\in I'$. If $\varphi_{|J}$ is constant, then $p_2'=\varphi(J)\subseteq I'$. In both cases, $p_2'\in I'$. We conclude that $p_2'$ is a common point  of the lines in the family $\{ I_j'\}$, so that $p_2'=p_1'$.

\end{proof}

\begin{thm}\label{7.10} Let $\mathbb{A}$ and $\mathbb{P}\,'\neq *$ be an affine space and a projective space over division rings $K$ and $K'$, respectively, with $\,|K|\geq 4\,$ or $\,|K|=3=\mathrm{char}\,K'$.

Any morphism of geometries $\varphi\colon\mathbb{A}\to\mathbb{P}\,'$, whose image is not contained in a line, extends to a unique partial morphism $\phi\colon\mathbb{P}\dashrightarrow\mathbb{P}\,'$, where $\mathbb{P}$ is the projective closure of $\mathbb{A}$.
\end{thm}
\begin{proof} Let $E\subset H$ be the set of exceptional points according to \ref{4.4b}. We define $\phi\colon\mathbb{P}-E\longrightarrow\mathbb{P}\,'$ by the formula
$$\phi(x):=\,\left\{\begin{aligned}
&\varphi(x)\qquad &\text{when }x\in\mathbb{A}\qquad\\
&\,x'\qquad &\text{when }x\in H\!-\!E
\end{aligned}\right.$$

Let us see that $\phi\colon\mathbb{P}\dashrightarrow\mathbb{P}\,'$ is a partial morphism of geometries, with exceptional subspace $E$. Since $\mathbb{P}$ is generated by lines, we have to prove
$$x_0\in x_1\vee x_2\ \Rightarrow\ \phi(x_0)\subseteq \phi(x_1)\vee\phi(x_2)$$
for any different points $x_0,x_1,x_2\in \mathbb{P}$. There are several cases according to the distribution of the points $x_0,x_1,x_2$ in $\mathbb{A}$ and $H$ :\smallskip

-- Case $x_0,x_1,x_2\in\mathbb{A}$. We apply that $\phi_{|\mathbb{A}}=\varphi$ is a morphism.\smallskip

-- Case $x_0,x_1\in\mathbb{A}$, $x_2\in H$. The condition $x_0\in x_1\vee x_2$ implies that $x_2$ is the  point at infinity of $x_0\vee x_1$. If $\varphi(x_0)\neq\varphi(x_1)$ then $x_2$ is not exceptional by Lemma \ref{4.3} and $\varphi(x_0),\varphi(x_1),x_2'$ are different collinear points by Lemma \ref{8.4}b; hence $\phi(x_0)=\varphi(x_0)\in\varphi(x_1)\vee x_2'=\phi(x_1)\vee\phi(x_2)$. If $\varphi(x_0)=\varphi(x_1)$, i.e., $\phi(x_0)=\phi(x_1)$, it clear that $\phi(x_0)\subseteq\phi(x_1)\vee\phi(x_2)$.\smallskip

-- Case $x_0,x_2\in\mathbb{A}$, $x_1\in H$. It similar to the previous one.\smallskip

-- Case $x_1,x_2\in\mathbb{A}$, $x_0\in H$. The condition $x_0\in x_1\vee x_2$ implies that $x_0$ is the  point at infinity of $x_1\vee x_2$. If $\varphi(x_1)\neq\varphi(x_2)$ then $x_0$ is not exceptional by Lemma \ref{4.3} and $\varphi(x_1),\varphi(x_2),x_0'$ are different collinear points by Lemma \ref{8.4}b; hence $\phi(x_0)=x_0'\in\varphi(x_1)\vee\varphi(x_2)=\phi(x_1)\vee\phi(x_2)$. If $\varphi(x_1)=\varphi(x_2)$ and $x_0$ is not excepcional then $\varphi(x_1)=\varphi(x_2)=x_0'$ by Lemma \ref{8.4}a, hence $\phi(x_1)=\phi(x_2)=\phi(x_0)$. If $x_0$ is exceptional then $\phi(x_0)=\emptyset\subseteq\phi(x_1)\vee\phi(x_2)$.\smallskip

-- Case $x_0,x_1,x_2\in H$. If $x_0\in E$ then $\phi(x_0)=\emptyset$ and it is clear. Let us assume $x_0\notin E$. If $x_1,x_2\notin E$ then $x_0\in x_1\vee x_2$ implies $\phi(x_0)\in\phi(x_1)\vee\phi(x_2)$ by Lemma \ref{8.5}. If $x_1\notin E$, $x_2\in E$, then $x_0\in x_1\vee x_2$ implies $(x_0\vee x_1)\cap E=x_2$, hence $x_0'=x_1'$ by Lemma \ref{6.6}, and $\phi(x_0)=x_0'=x_1'=\phi(x_1)\in\phi(x_1)\vee\phi(x_2)$. The case $x_2\notin E$, $x_1\in E$ is similar. If $x_1,x_2\in E$ then $x_0\in x_1\vee x_2$ implies $x_0\in E$ by Proposition \ref{4.4a}, against the assumption on $x_0$.

\end{proof}

Remark that the above partial morphism $\phi\colon\mathbb{P}\dashrightarrow\mathbb{P}'$ is a partial projective morphism by the Fundamental Theorem \ref{5.6}. A related theorem is obtained in (\cite{Frank} Satz 2)  replacing $\mathbb{A}$ by a non empty open set of $\mathbb{P}$ with respect to a linear topology. \goodbreak

\begin{cor} With the hypotheses of the theorem, if $\varphi\colon\mathbb{A}\to\mathbb{P}\,'$ is injective then the extension  $\phi\colon\mathbb{P}\to\mathbb{P}\,'$ is globally defined.
\end{cor}
\begin{proof} We have $E=\emptyset$: If there is an exceptional point $p$, then for any line $L\subset\mathbb{A}$ with the point $p$ at infinity we have that $\varphi_{|L}$ is constant, so contradicting that $\varphi$ is injective.

\end{proof}

\begin{cor}\label{5.11} With the hypotheses of the theorem, if $\varphi\colon\mathbb{A}\to\mathbb{P}\,'$ is an embedding then the extension $\phi\colon\mathbb{P}\to\mathbb{P}\,'$ so is an embedding.
\end{cor}
\begin{proof} If $\text{dim}\,\mathbb{A}=n<\infty$, then we have $\text{dim}\,\phi(\mathbb{P})\leq\text{dim}\,\mathbb{P}=n$ and $\text{dim}\,\phi(\mathbb{P})\geq \text{dim}\,\varphi(\mathbb{A})=\text{dim}\,\mathbb{A}=n$, so that $\text{dim}\,\phi(\mathbb{P})=n$. By Proposition \ref{1.5}b, $\phi\colon\mathbb{P}\to\phi(\mathbb{P})$ is an isomorphism, that is to say, $\phi\colon\mathbb{P}\to\mathbb{P}\,'$ is an embedding.

In the general case we put $\mathbb{A}$ as the union of all affine subspaces of finite dimension: $\mathbb{A}=\bigcup\mathbb{A}_i$. Then $\mathbb{P}=\bigcup\mathbb{P}_i$, where $\mathbb{P}_i$ is the projective closure of $\mathbb{A}_i$. By the former case, $\phi\colon\mathbb{P}_i\to\phi(\mathbb{P}_i)$ is an isomorphism for any index $i$, hence $\,\phi\colon\mathbb{P}=\bigcup\mathbb{P}_i\longrightarrow\bigcup\phi(\mathbb{P}_i)=\phi(\mathbb{P})\,$ is an isomorphism.

\end{proof}

This result was obtained in \cite{Rigby} in the case $\text{dim}\,\mathbb{A}=2$ and in \cite{BR} in the case $\text{dim}\,\mathbb{A}<\infty$. It may be generalized replacing $\mathbb{A}$ by an open set of $\mathbb{P}$ with respect to a linear topology (see \cite{Frank} Satz 1). For the existence of embeddings between affine or projective spaces, over commutative fields, see \cite{T1}-\cite{T4}.

\begin{nada}{\it Counterexamples}\smallskip

-- When $|K|=3$, it is well known that the affine plane $\mathbb{A}_2$ may be embedded in the complex projective plane $\mathbb{P}_2'$ as the set of nine inflection points 
of a non singular cubic, since any line passing through two inflection points also passes through a third one (this property is related with the group law of the cubic). 
This embedding may be not extended to the projective closure of $\mathbb{A}_2$ because there are no homomorphisms $K\to\mathbb{C}$.\smallskip

-- Let $\mathbb{A}_3$ be the 3-dimensional affine space over $K=\mathbb{Z}_2$; this space has eight points $p_1,\ldots ,p_8$ (vertices of a cube) and each line has two points. Let $\mathbb{P}'_2$ be the projective plane over a field $K'$ with $|K'|\ge 7$; so that the non singular conic $C\subset \mathbb{P}_2'$ given by the homogeneous equation $x_1^2-x_2x_3=0$ has $\ge 8$ points. Let us consider an octagon $X=\{ p'_1,\ldots ,p'_8\} \subseteq C$ inscribed in the conic. Since there are not three collinear vertices in the octagon, it is easy to check that the injective map $\mathbb{A}_3\to \mathbb{P}'_2$, $p_i\mapsto p'_i$, is a morphism of geometries. If we order the points of $\mathbb{A}_3$ so that the lines $p_1p_2$, $p_3p_4$ and $p_5p_6$ are parallel, and we take an octagon such that the sides $p'_1p'_2$, $p'_3p'_4$, $p'_5p'_6$ are not concurrent, then it is clear that the morphism $\mathbb{A}_3\to \mathbb{P}'_2$ does not extend to the projective closure of $\mathbb{A}_3$.
\end{nada}

\section{Morphisms $\mathbb{A}\to\mathbb{A}'$}

Now we use Theorem \ref{7.10} to prove that the morphisms of geometries $\varphi\colon\mathbb{A}\to\mathbb{A}'$ are fractional semiaffine maps.

\begin{thm}\label{6.1} Let $\mathbb{A}$ and $\mathbb{A}'\neq *$ be affine spaces over division rings $K$ and $K'$, respectively, with $|K|\geq 4$ or $|K|=3=\mathrm{char}\, K'$.

Let $\varphi\colon\mathbb{A}\to\mathbb{A}'$ be a map such that the image is not contained in a line. The following conditions are equivalent:\smallskip

a) $\,\varphi\colon\mathbb{A}\to\mathbb{A}'$ is a morphism of geometries: for all $p_0,p_1,p_2\in\mathbb{A}$ it holds $\,p_0\in p_1\vee p_2\quad\Rightarrow\quad \varphi(p_0)\in\varphi(p_1)\vee\varphi(p_2)$.\smallskip

b) $\,\varphi\colon\mathbb{A}\to\mathbb{A}'\subset\mathbb{E}'\,$ is a fractional semiaffine morphism: $\varphi=\frac{1}{\omega}\Phi\,$ for some semiaffine morphisms $\Phi\colon\mathbb{A}\to\mathbb{E}'$, $\omega\colon\mathbb{A}\to K'$.
\end{thm}
\begin{proof} ($a\Rightarrow b$). By theorems \ref{7.10} and \ref{5.6}, the morphism $\varphi\colon\mathbb{A}\to\mathbb{A}'\subset\mathbb{P}\,'$ extends to a partial projective morphism  $\phi\colon\mathbb{P}=\mathbb{P}(\mathbb{E})\,\dashrightarrow\,\mathbb{P}(\mathbb{E}')=\mathbb{P}\,'$, $$\xymatrixcolsep{1pc}\xymatrixrowsep{2pc}
\xymatrix{
\mathbb{A} \ar[d]_{\varphi} &  \subset & \mathbb{P}(\mathbb{E}) \ar@{-->}[d]^{\phi}\\
\mathbb{A}' &  \subset & \mathbb{P}(\mathbb{E}')
}$$
The morphism $\phi$ is defined by a semilinear map $\Phi\colon \mathbb{E}\to\mathbb{E}'$. Now, the above commutative square shows that the following diagram
$$\xymatrixcolsep{1pc}\xymatrixrowsep{2pc}
\xymatrix{
\mathbb{A} \ar[d]_{\varphi} &  \subset & \mathbb{E} \ar[d]^{\Phi}\\
\mathbb{A}' &  \subset & \mathbb{E}'
}$$
is commutative modulo a scalar factor:
\begin{equation}\label{eq4}
\Phi(x)\,=\,\lambda'_x\cdot\varphi(x)\qquad\text{with}\quad 0\neq \lambda_x'\in K'
\end{equation}
for all $x\in\mathbb{A}$.

Let $\omega'\colon\mathbb{E}'\to K'$ be the unique $K'$-linear map such that $\omega'_{|\mathbb{A}'}=1$.  We define $\,\omega:=\omega'\circ\Phi_{|\mathbb{A}}:\mathbb{A}\to K'$; since $\omega'$ is linear, 
the semiaffine morphisms $\Phi,\omega$ have the same associated homomorphism $K\to K'$.

For all $x\in\mathbb{A}$ we have
$$\omega(x)\,=\,\omega'(\Phi(x))\,\overset{(\ref{eq4})}{=}\,\lambda'_x\,\omega'(\varphi(x))\,\overset{\omega'_{|\mathbb{A}'}=1}{=}\,\lambda'_x\cdot 1\,\neq\, 0$$
so that (\ref{eq4}) states that
$$
\Phi(x)\,=\,\omega(x)\cdot\varphi(x)\qquad\quad\forall\ x\in\mathbb{A}$$
hence $\varphi=\frac{1}{\omega}\Phi$ is a fractional semiaffine morphism.\smallskip

($b\Rightarrow a$). It is Proposition \ref{4.15}.

\end{proof}

This theorem was obtained by W. Zick in an unpublished paper \cite{Zick2}. According to it, the result was previously proved by G. Martin when  $K$ and $K'$ are commutative.\medskip

Combining with \ref{5.1a} we obtain the following generalization of the classical Fundamental Theorem of real affine geometry,

\begin{cor} Let $\varphi\colon\mathbb{A}\to\mathbb{A}'$ be a map between real affine spaces, such that the image is not contained in a line. Then $\varphi$ is an affine morphism if and only if $\varphi$ is a morphism of geometries.
\end{cor}

\noindent{\it Comments}\smallskip

-- The above corollary may be easily obtained from a result of Lenz \cite{Lenz}: {\it Let  $\mathbb{P}$, $\mathbb{P}'$ be real projective spaces. Any morphism of geometries $\varphi\colon U\to\mathbb{P}\,'$, defined on an open subset $U\subseteq\mathbb{P}$, such that the image is not contained in a line, extends to a unique partial projective morphism $\phi\colon\mathbb{P}\dashrightarrow\mathbb{P}\,'$.}\medskip

-- The isomorphisms of geometries $\varphi\colon \mathbb{A}\to\mathbb{A}'$ are the collineations. By the Fundamental Theorem of Affine Geometry, the collineations $\varphi\colon \mathbb{A}\to\mathbb{A}'$ are just the semiaffinities (assuming $\text{dim}\,\mathbb{A}\geq 2$ and $|K|\neq 2$). This classical result may be obtained easily as a consequence of Theorem \ref{6.1}. An extension of the fundamental theorem, which geometrically characterizes semiaffine morphisms, was obtained by Zick \cite{Zick1}, see also \cite{Sancho1}.\smallskip

-- A map $\varphi\colon \mathbb{A}\to\mathbb{A}'$ is called a {\it lineation} if the images
by $\varphi$ of any three collinear points are collinear. It is a weaker notion than morphism of geometries and it has been studied in the literature, specially in dimension $2$ (see \cite{CV}). A version of \ref{5.11} for injective lineations (and finite dimension) is given in \cite{BR}. See also \cite{ChP} for a version of the Fundamental Theorem using surjective lineations.



\end{document}